\documentstyle[12pt]{article}

\newfont{\frakturfont}{eufm10 scaled\magstep1}
\newfont{\smallfrakturfont}{eufm7 scaled\magstep1}
\newfont{\bbbfont}{msbm10 scaled\magstep1}
\newfont{\footbbbfont}{msbm8 scaled\magstep1}
\newfont{\smallbbbfont}{msbm7 scaled\magstep1}
\newfont{\tinybbbfont}{msbm5 scaled\magstep1}

\advance\textheight by -\headheight
\advance\textheight by -\headsep
\makeatletter

\long\def\@makecaption#1#2{
 \vskip 10pt
 \setbox\@tempboxa\hbox{#1. #2}
 \ifdim \wd\@tempboxa >\hsize #1. #2\par \else \hbox
to\hsize{\hfil\box\@tempboxa\hfil}
 \fi}

\def\thebibliography#1{\section*{References\@mkboth
 {REFERENCES}{REFERENCES}}\small\list
 {\arabic{enumi}.}{\settowidth\labelwidth{[#1]}\leftmargin\labelwidth
 \advance\leftmargin\labelsep
 \usecounter{enumi}}
 \def\newblock{\hskip .11em plus .33em minus .07em}
 \sloppy\clubpenalty4000\widowpenalty4000
 \sfcode`\.=1000\relax}

\newcommand{\text}[1]{\mathchoice{\mbox{\rm #1}}{\mbox{\rm #1}}
    {\mbox{\scriptsize\rm #1}}{\mbox{\tiny\rm #1}}}
\newcommand{\operatorname}[1]{\mathop{\rm #1}\nolimits}

\newcommand{\pfit}[1]{{\it #1:\/}}
\newcommand{\qedsymbol}{Q.E.D.}
\newcommand{\eqref}[1]{(\ref{#1})}

\newenvironment{pf*}{\noindent\pfit }{\qedsymbol\par\par\medskip\par}

\newenvironment{pmatrix}{\left(\begin{array}{ccccc}}{\end{array}\right)}

\def\maketitle{\par
 \begingroup
 \def\thefootnote{\fnsymbol{footnote}}
 \def\@makefnmark{\hbox
 to 0pt{$^{\@thefnmark}$\hss}}
 \if@twocolumn
 \twocolumn[\@maketitle]
 \else \newpage
 \global\@topnum\z@ \@maketitle \fi\thispagestyle{empty}\@thanks
 \endgroup
 \setcounter{footnote}{0}
 \let\maketitle\relax
 \let\@maketitle\relax
 \gdef\@thanks{}\gdef\@author{}\gdef\@title{}\let\thanks\relax}

\def\@address{}
\def\address#1{\gdef\@address{#1}}
\def\@email{}
\def\email#1{\gdef\@email{#1}}
\def\at{@}
\def\trailer{\bigskip\par\noindent \footnotesize {\sc \@address}
   \par\noindent {\it E-mail address:} \@email}
\date{}
\def\dedicatory#1{\def\@date{\normalsize\it#1}}
\def\subjclass#1{\def\@thefnmark{}\@footnotetext{1991
    {\it Mathematics Subject Classification.} #1}}
\def\keywords#1{\def\@thefnmark{}\@footnotetext{
    {\it Key words and phrases.} #1}}
\def\thanks#1{\def\@thefnmark{}\@footnotetext{#1}}

\def\@begintheorem#1#2{\trivlist \item[\hskip \labelsep{\sc #1\ #2}]\it}
\def\@opargbegintheorem#1#2#3{\trivlist
      \item[\hskip \labelsep{\sc #1\ #2\ (#3)}]\it}

\makeatother

\newtheorem{definition}{Definition}

\newtheorem{ConeConjecture}{The Cone Conjecture}
  
\newtheorem{ConvergenceConjecture}{The Convergence Conjecture}
  
\newtheorem{mirrorconjecture}{The Mathematical Mirror Symmetry Conjecture}

\newtheorem{minconjecture}{The Minimal Compactification Conjecture}

\newcommand{\frak}[1]{\mathchoice{\mbox{\frakturfont #1}}{\mbox{\frakturfont
#1}}
    {\mbox{\smallfrakturfont #1}}{\mbox{\smallfrakturfont #1}}}

\newcommand{\Bbb}[1]{\mathchoice{\mbox{\bbbfont #1}}{\mbox{\bbbfont #1}}
    {\mbox{\smallbbbfont #1}}{\mbox{\tinybbbfont #1}}}

\newcommand{\C}{\Bbb C}
\newcommand{\E}{\Bbb E}
\renewcommand{\P}{\Bbb P}
\newcommand{\Q}{\Bbb Q}
\newcommand{\R}{\Bbb R}

\newcommand{\Z}{\Bbb Z}
\newcommand{\cB}{{\cal B}}
\newcommand{\cC}{{\cal C}}
\newcommand{\cD}{{\cal D}}
\newcommand{\cF}{{\cal F}}
\newcommand{\cK}{{\cal K}}
\newcommand{\cM}{{\cal M}}
\newcommand{\cO}{{\cal O}}
\newcommand{\cP}{{\cal P}}
\newcommand{\cS}{{\cal S}}

\newcommand{\cY}{{\cal Y}}
\renewcommand{\aa}{\frak A}
\newcommand{\D}{\Delta} 
\newcommand{\h}{\frak H} 
\renewcommand{\O}{\frak O}

\newcommand{\rtimes}{\mathbin{\mbox{\bbbfont\char"6F}}}
\newcommand{\semidirect}{\rtimes}
\newcommand{\dlog}{d\mskip0.5mu\log}
\newcommand{\flt}{\nabla}
\newcommand{\stminus}{-}

\newcommand{\suchthat}{ :  }

\newcommand{\Aff}{\operatorname{Aff}}
\newcommand{\Aut}{\operatorname{Aut}}
\newcommand{\DLP}{\operatorname{DLP}}
\newcommand{\GL}{\operatorname{GL}}
\newcommand{\PGLplus}{\operatorname{PGL\mbox{\mathsurround=0pt$^+$}}}
\newcommand{\PSL}{\operatorname{PSL}}
\newcommand{\SL}{\operatorname{SL}}
\newcommand{\SU}{\operatorname{SU}}
\renewcommand{\Im}{\operatorname{Im}}
\newcommand{\Ker}{\operatorname{Ker}}
\newcommand{\Lie}{\operatorname{Lie}}
\newcommand{\rk}{\operatorname{rk}}
\newcommand{\diag}{\operatorname{diag}}

\newcommand{\lrp}{locally rational polyhedral}
\newcommand{\rp}{rational polyhedral}
\newcommand{\SBB}{Satake-Baily-Borel}
\newcommand{\NLSM}{nonlinear sigma-model}

\begin{document}

\title{Compactifications of moduli spaces \\inspired by mirror symmetry}
\author{David R. Morrison}
\date{}

\maketitle

\address{Department of Mathematics, Duke University,
Durham, NC 27708-0320 USA, and \\School of Mathematics, Institute for
Advanced Study, Princeton, NJ 08540 USA.}
\email{drm\at math.duke.edu}

The study of moduli spaces by means of
the period mapping has found its greatest
success for moduli spaces of varieties with trivial canonical
bundle, or more generally, varieties with Kodaira dimension
zero.  Now these moduli spaces play a pivotal
r\^ole in the classification theory of algebraic
varieties, since varieties with nonnegative Kodaira dimension
which are not of general type admit birational fibrations by
varieties of Kodaira dimension zero.  Since such fibrations typically
include singular fibers as well as smooth ones,  it is
important to understand how to compactify the corresponding moduli spaces
(and if possible, to give geometric interpretations to the
boundary of the compactification).  Note that because of the
possibility of blowing up along the boundary, abstract compactifications
of moduli spaces are far from unique.

The hope that the period mapping could be used to construct
compactifications of moduli spaces was given concrete expression
in some conjectures of Griffiths \cite[\S 9]{gr-bull} and others in the
late 1960's.  In particular, Griffiths conjectured that there would
be an analogue of the \SBB\ compactifications of arithmetic quotients of
bounded symmetric domains, with some kind of ``minimality'' property
among compactifications.  Although there has been much progress since
\cite{gr-bull} in
understanding the behavior of period mappings near the boundary
of moduli,
 compactifications of this type have not been constructed, other than
in special cases.

In the case of algebraic
K3 surfaces, the moduli spaces themselves are arithmetic
quotients of bounded symmetric domains,  so each has a minimal (\SBB)
compactification.  In studying the moduli spaces for K3 surfaces of
low degree in the early 1980's, Looijenga \cite{Looijenga} found
that the \SBB\ compactification needed to be blown up slightly in order
to give a good geometric interpretation to the boundary.
He introduced a class of compactifications, the {\em semi-toric
compactifications}, which includes the ones with a good geometric
interpretation.

In higher dimension, the moduli spaces are not expected to be arithmetic
quotients of symmetric domains, so different techniques are needed.
The study of these moduli spaces has received renewed attention
recently, due to the discovery by theoretical physicists of
a phenomenon called ``mirror symmetry''.  One of the predictions of
mirror symmetry is
that the moduli space for a variety with trivial canonical bundle,
which parameterizes the possible complex structures on the underlying
differentiable manifold,
should {\em also}\/ serve as the parameter space for a very different kind of
structure on a ``mirror partner''---another
  variety with trivial canonical
bundle.  This alternate description of the moduli space turns out to be
well-adapted to analysis by Looijenga's techniques; we carry out that
analysis here.

In the physicists' formulation, one fixes a differentiable manifold
$X$  which admits complex structures with
trivial canonical bundle (a ``Calabi-Yau manifold''),
and studies something called  \NLSM s on $X$.
Such an object can be determined by specifying both a complex structure
on $X$, and some ``extra structure'' (cf.~\cite{guide});
the moduli space of interest
to the physicists parameterizes the choice of both.
The r\^oles of the ``complex structure'' and ``extra structure'' subspaces
of this parameter space are reversed when $X$ is replaced by a mirror
partner.

Most aspects of mirror symmetry must be regarded as conjectural by
mathematicians at the moment, and in this paper we conjecture much
more than we prove.
In a companion paper \cite{htams},
we consider formally degenerating variations of Hodge
structure near normal crossing boundary points of the
moduli space, and describe a conjectural link to the
numbers of rational curves of various degrees on a mirror partner.
In the present paper, we extend these
considerations to boundary points which are {\em not}\/
of normal crossing type,
and formulate a mathematical mirror symmetry conjecture in greater
generality.
  In addition, we find that when studied from the mirror perspective,  a
``minimal'' partial compactification of the moduli space---analogous
to the \SBB\ compactification---appears very natural,
provided that several conjectures about the mirror partner hold.

One of our conjectures is a simple and compelling statement about
the K\"ahler cone of Calabi-Yau varieties.  If true, it clarifies
the r\^ole of some of the ``infinite discrete'' structures on
such a variety, which nevertheless seem to be finite modulo automorphisms.
We have verified this conjecture in a nontrivial case in joint work
with A.~Grassi \cite{grassi-morrison}.

The plan of the paper is as follows.
In the first several sections, we review Looijenga's compactifications,
describe a concrete example, and add a refinement to the theory in
the form of a flat connection on the holomorphic cotangent bundle of
the moduli space.  We then turn to the description of the larger
moduli spaces of interest to physicists, and analyze certain boundary
points of those spaces.  Towards the end of the paper, we explore the
mathematical implications of mirror symmetry in constructing
 compactifications of moduli
spaces.
We close by discussing some evidence for mirror symmetry which
(in hindsight) was available in 1979.

\section{Semi-toric compactifications} \label{STC}

The first methods for compactifying
 arithmetic quotients of bounded symmetric domains were
found
by Satake \cite{satake1} and Baily-Borel \cite{baily-borel}.
The compactification  produced by their methods,
often called the {\em \SBB\ compactification},
adds a ``minimal'' amount to the quotient space
in completing it to a compact complex analytic space.
This minimality can be made quite precise, thanks to the
Borel extension theorem \cite{borel-extension}
which guarantees that for a given
quotient of a bounded symmetric domain by an arithmetic group,
any compactification whose boundary is a divisor
with normal crossings
 will map to the \SBB\ compactification
(provided that the arithmetic group is torsion-free).

\SBB\ compactifications have rather bad singularities on their
boundaries,  so they are difficult to study in detail.
Explicit resolutions
of singularities for these compactifications
were constructed in special cases by
Igusa \cite{igusa},
 Hemperly \cite{hemperly}, and Hirzebruch \cite{hirzebruch};
 the general case was subsequently treated
by Satake \cite{satake2} and Ash et.~al \cite{AMRT}.
The methods of \cite{AMRT} produce what are usually called
{\em Mumford compactifications}---these are smooth, and have a
divisor with normal crossings on the
boundary, but unfortunately many choices must be made
in their construction.  The \SBB\ compactification, on the other
hand, is canonical.

Some years later, Looijenga \cite{Looijenga}
generalized both the \SBB\ and the
Mumford compactifications by means of a construction which can be
applied widely, not just in the case of arithmetic quotients of
bounded symmetric domains.
Looijenga's construction gives partial compactifications of certain
quotients of tube domains by discrete group actions.  A {\em tube domain}\/
is the set of points in a complex vector space whose imaginary parts
are constrained to lie in a specified cone.  Whereas Ash et
al.\ \cite{AMRT} had only considered homogeneous self-adjoint cones,
Looijenga showed that analogous constructions could be made in a
more general context.

The starting point is a free $\Z$-module $L$ of finite rank,
and the real vector space $L_{\R}:=L\otimes\R$ which it spans.
A convex cone $\sigma$ in $L_{\R}$ is {\em strongly convex}\/
if $\sigma\cap(-\sigma)\subset\{0\}$.  A convex cone is {\em generated}\/
by the set $S$ if every element in the cone can be written as
a linear combination of the elements of $S$ with nonnegative coefficients.
And a convex cone is {\em rational polyhedral}\/
 if it is generated by a finite
subset of the rational vector space $L_{\Q}:=L\otimes\Q$.

Let $\cC\subset L_{\R}$
be an open strongly convex cone, and
let $\Gamma\subset\Aff(L)$ be a group of affine-linear transformations of $L$
which contains the translation subgroup $L$ of $\Aff(L)$.
If the linear part $\Gamma_0:=\Gamma/L\subset\GL(L)$
of $\Gamma$ preserves the cone $\cC$, then the group $\Gamma$ acts
on the tube domain $\cD:=L_{\R}+i\,\cC$.
We wish to partially compactify the quotient space $\cD/\Gamma$,
including limit points for
all paths moving out towards infinity in the tube domain.

Looijenga formulated a condition which guarantees the existence of
partial compactifications of this kind.  Let
 $\cC_+$ be the convex hull of $\overline{\cC}\cap L_{\Q}$.
Following
\cite{Looijenga},
we say that $(L_{\Q},\cC,\Gamma_0)$ is {\em admissible}\/
if there exists a rational polyhedral cone
$\Pi\subset\cC_+$ such that $\Gamma_0.\Pi=\cC_+$.
Given an admissible triple $(L_{\Q},\cC,\Gamma_0)$,
the (somewhat cumbersome)
 data needed to specify one of Looijenga's partial compactifications is as
follows.\footnote{We have modified Looijenga's definition
slightly, so that the use of the term  ``face'' is the standard
one (cf.~\cite{rockafellar}):  a subset
$\cF$  of a convex set $\cS$ is a {\em face}\/ of $\cS$
if every closed line segment
in $\cS$ which has one of its relative interior points lying in $\cF$
also has both endpoints lying in $\cF$.}

\begin{definition} {\rm \cite{Looijenga}}
  A\/ {\em \lrp\ decomposition} of $\cC_+$ is a collection
$\cP$ of strongly convex cones
such that
\begin{itemize}
\item[(i)]  $\cC_+$ is the disjoint union
of the cones belonging to $\cP$,

\item[(ii)] for every $\sigma\in\cP$, the $\R$-span of $\sigma$ is
defined over $\Q$,

\item[(iii)] if $\sigma\in\cP$, if $\tau$ is the relative interior of
a nonempty
face of the closure of $\sigma$, and if $\tau\subset\cC_+$, then
$\tau\in\cP$, and

\item[(iv)] if $\Pi$ is a \rp\ cone in $\cC_+$, then $\Pi$ meets only
finitely many members of $\cP$.

\end{itemize}
\end{definition}

(The decomposition $\cP$ is called {\em \rp}\/ if all the cones in $\cP$ are
relative interiors of \rp\ cones.  This is
the same notion which appears in toric geometry \cite{Fulton,Oda},
except that the cones appearing in $\cP$ as formulated here
are the relative interiors
of the cones appearing in that theory.)

For each $\Gamma_0$-invariant
\lrp\ decomposition $\cP$ of $\cC_+$, there is a partial compactification
 of $\cD/\Gamma$ called the
{\em semi-toric (partial) compactification associated
to $\cP$}.
This partial compactification has the form $\widehat{\cD}(\cP)/\Gamma$, where
$\widehat{\cD}(\cP)$ is the disjoint union of certain strata $\cD(\sigma)$
associated to the cones $\sigma$ in the decomposition.
The complex dimension of the stratum
$\cD(\sigma)$ coincides with the real codimension
of the cone $\sigma$ in $L_{\R}$; in particular, the open cones in $\cP$
correspond to the $0$-dimensional strata in $\widehat{\cD}(\cP)$.
The delicate
points in the construction are the specification
of a topology on $\widehat{\cD}(\cP)$, and the proof that
the quotient space $\widehat{\cD}(\cP)/\Gamma$ has a natural structure of
a normal complex analytic space.
For more details, we refer the reader to \cite{Looijenga} or \cite{Sterk}.

The construction has the property that
if $\cP'$ is a refinement of $\cP$, then there is a dominant morphism
$\widehat{\cD}(\cP')/\Gamma\to\widehat{\cD}(\cP)/\Gamma$.
Blowups of the boundary can be realized in this way.

A bit more generally, we can partially compactify finite
covers $\cD/\Gamma'$ of $\cD/\Gamma$,
built from $L'\subset L$ of finite index,
$\Gamma_0'\subset\GL(L')\cap\Gamma_0$ of finite index in $\Gamma_0$,
and
$\Gamma':=L'\semidirect\Gamma_0'$,
by specifying  a $\Gamma_0'$-invariant \lrp\ decomposition $\cP'$ of $\cC_+$.

There are two extreme cases of a semi-toric compactification.
The {\em \SBB\ decomposition}\/ $\cP_{\text{SBB}}$ consists of all
relative interiors of nonempty faces of $\cC_+$.  The resulting
(partial) compactification $\widehat{\cD}(\cP_{\text{SBB}})/\Gamma$
is the {\em \SBB-type compactification of $\cD/\Gamma$}.  This is
``minimal'' among semi-toric compactifications in an obvious combinatorial
sense; I do not know whether a more precise
 analogue of the Borel extension theorem
holds in this context.
The strata added to $\cD/\Gamma$ include a unique $0$-dimensional
stratum $\cD(\cC)$, which serves as a distinguished boundary point.

At the other extreme, if
every cone $\sigma\in\cP$ is the relative interior of a \rp\ cone
$\overline{\sigma}$ which is
generated by a subset of a
 basis of $L$, then the associated partial compactification is
smooth, and the compactifying set is a divisor with normal crossings.
We call this a
{\em Mumford-type}\/ semi-toric compactification.
We will spell out the structure of the compactification more
explicitly in this case, giving an alternative description of
$\widehat{\cD}(\cP)/\Gamma$.

We can think of producing a Mumford-type semi-toric compactification
in two steps.  In the first step, we construct a partial compactification
 $\widehat{\cD}(\cP)/L$ of $\cD/L$ which is $\Gamma_0$-equivariant;
in the second step we
recover $\widehat{\cD}(\cP)/\Gamma$ as the quotient of
$\widehat{\cD}(\cP)/L$ by $\Gamma_0$.

The first step is done one cone at a time.
Given $\sigma\in\cP$, there is a
basis $\ell^1,\dots,\ell^r$ of $L$ such that
\[\sigma=\R_{>0}\ell^1+\cdots+\R_{>0}\ell^k\quad\text{for some $k\le r$}.\]
Let $\{z_j\}$ be complex coordinates dual to $\{\ell^j\}$, so that
$z=\sum z_j\ell^j$ represents a general element of $L_{\C}$.
Consider the set $\cD_{\sigma}:=L_{\R}+i\,\sigma$.  Translations by
the lattice $L$ preserve $\cD_{\sigma}$, and coordinates on the
quotient $\cD_{\sigma}/L\subset L_{\C}/L$ can be given by
 $w_j=\exp(2\pi i\,z_j)$.  In terms of those coordinates, $\cD_{\sigma}/L$
can be described as
\[\cD_{\sigma}/L=\{w\in\C^r\suchthat  0<|w_j|<1 \text{ for $j\le k$},
|w_j|=1 \text{ for $j>k$}\}.\]
We partially compactify this to
\[(\cD_{\sigma}/L)^-:=\{w\in\C^r\suchthat  0\le|w_j|<1 \text{ for $j\le k$},
|w_j|=1 \text{ for $j>k$}\}.\]
(We have suppressed the $\sigma$-dependence of $\ell^j$, $z_j$, $w_j$
to avoid cluttering up the notation.)
We call any $w\in(\cD_{\sigma}/L)^-$ with $w_j=0$ for $j\le k$
a {\em distinguished
limit point}\/ of $\cD_{\sigma}/L$.
Note that any path in $\cD_{\sigma}$
along which $\Im(z_j)\to\infty$ for all
$j\le k$, maps to a path in $\cD_{\sigma}/L$
which approaches such a distinguished limit point.
The set $\DLP(\sigma)$ of distinguished limit points is a subset of the
stratum $\widehat{\cD}(\sigma)$, and is a compact real torus of dimension
$\dim_{\R}\DLP(\sigma)=r-k=\dim_{\C}\widehat{\cD}(\sigma)$.
When $k=r$, the distinguished limit point is unique, and it coincides
with the $0$-dimensional stratum $\cD(\sigma)$ of $\widehat{\cD}(\cP)$.

The partial compactification $\widehat{\cD}(\cP)/L$ can now be described
as a disjoint union of the $(\cD_{\sigma}/L)^-$'s, with
$(\cD_{\tau}/L)^-$ lying in the closure of $(\cD_{\sigma}/L)^-$
whenever $\tau$ is the relative interior of
a face of $\overline{\sigma}$.  This space
$\widehat{\cD}(\cP)/L$ is smooth and simply-connected,
and the induced action of $\Gamma_0$
 on it has no fixed points.
The action of $\Gamma_0$
permutes the various $(\cD_{\sigma}/L)^-$'s, a finite number of which
serve to cover $\widehat{\cD}(\cP)/\Gamma$ after we take the
quotient by $\Gamma_0$.

The structure of $(\cD_\sigma/L)^-$ near the distinguished limit point
when $k=r$ can be formalized in the following way.  For  a
complex manifold $T$,  we
 say that $p$ is a {\em maximal-depth normal crossing point
of $B\subset T$}\/
if there is an open neighborhood $U$ of $p$ in $T$
and an isomorphism $\varphi:U\to\D^r$
such that $\varphi(U\cap(T{-}B))=(\D^*)^r$ and $\varphi(p)=(0,\dots,0)$,
where $\D$ is the unit disk, and $\D^*:=\D{\stminus}\{0\}$.
There are thus $r$ local components $B_j:=\varphi^{-1}(\{v_j=0\})$ of
$B\cap U$, with $p=B_1\cap\cdots\cap B_r$, where $v_j$ is a coordinate
on the $j^{\text{th}}$ disk.

\section{Cusps of Hilbert modular surfaces}\label{HMS}

We now give an example to illustrate the
construction in the previous section:
the cusps of Hilbert modular surfaces, as
analyzed by Hirzebruch \cite{hirzebruch} and by
Mumford in the first chapter of \cite{AMRT}.
Let $\PGLplus(2,\R)=\PSL(2,\R)$ act
by fractional
linear transformations on the upper half plane $\h$.
Let $K$ be a real quadratic field with ring of integers $\O_K$,
and let $\PGLplus(2,K)$ be the group of invertible
$2\times2$ matrices with entries in $K$
whose determinant is mapped to a positive number
under  both embeddings
of $K$ into $\R$, modulo scalar multiples of the identity matrix.
The map
 $\Phi:K\to\R^2$  given by the two embeddings of $K$ into $\R$
induces an action of $\PGLplus(2,K)$ on $\h\times\h$.

A {\em Hilbert modular surface}\/
is an algebraic surface of
 the form $\h\times\h/\Gamma$ for some  arithmetic group
$\Gamma\subset\PGLplus(2,K)$
(that is, a group commensurable with $\PGLplus(2,\O_K)$),
often assumed to be torsion-free.
 The \SBB\ compactification
of  a Hilbert modular surface adds a finite number of compactification points,
called {\em cusps}.  Small deleted neighborhoods of such points
have inverse images in $\h\times\h$ whose $\Gamma$-stabilizer is
a parabolic subgroup $\Gamma_{\text{par}}$ of the form
\[\Gamma_{\text{par}}=\{\begin{pmatrix}\varepsilon^k & a \\ 0 & 1
\end{pmatrix}
 \suchthat k\in\Z,a\in\aa
\},\]
where $\aa\subset\O_K$ is an ideal, and $\varepsilon\in\O_K^\times$
is a totally positive unit such that $\varepsilon\,\aa=\aa$.  We
can analyze a neighborhood of a cusp by studying appropriate
partial compactifications of $\h\times\h/\Gamma_{\text{par}}$.

\setlength{\unitlength}{1 true cm}

\begin{figure}[t]
\begin{center}
\begin{picture}(10,5)(0,0)
\put(2,3){\mbox{Insert here the figure from p.~52 of \cite{AMRT}}}
\end{picture}
\end{center}
\caption{}
\end{figure}

The elements in $\Gamma_{\text{par}}$ with $k=0$ form the  translation
subgroup, which we identify with $\aa$.  This
is a free abelian group of rank 2. Let $(\alpha,\alpha')$,
$(\beta,\beta')$ be a $\Z$-basis of $\Phi(\aa)$.
Define a map $\h\times\h\to\C^2$
by
\[(w_1,w_2)\mapsto\frac1{\alpha\beta'-\alpha'\beta}
(\beta'w_1-\beta w_2,-\alpha'w_1+\alpha w_2),\]
and let $\cD$ denote the image of $\h\times\h$ in $\C^2$.
Under this map,
$\Phi(\aa)$ is sent to the standard lattice $L:=\Z^2$,
and $\Phi(\Gamma_{\text{par}})$ is sent to a subgroup of
$\Aff(L)$ with the translation subgroup $\aa$ of $\Gamma_{\text{par}}$
mapped to the translation subgroup $L$ of $\Aff(L)$.
As in section \ref{STC}, we form the quotient in two steps:
first take the quotient $\h\times\h/\aa=\cD/L$,
and then take the quotient of
the resulting space by the group $\Gamma_0:=\Gamma_{\text{par}}/\aa$.

Mumford shows how to partially compactify the space
$\cD/L\subset L\otimes\C^*=(\C^*)^2$ in a $\Gamma_0$-equivariant
way, so that the quotient by $\Gamma_0$ gives the desired
 partial compactification of $\h\times\h/\Gamma_{\text{par}}$.
The map of $\h\times\h\to\C^2$ was designed so that the image
would be a tube domain
$\cD:=\R^2+i\,\cC$,
where $\cC$ is the cone
\[\cC=\{(y_1,y_2)\suchthat  \alpha y_1 + \beta y_2>0,
\alpha'y_1+\beta'y_2>0\}.\]
The boundary lines of the closure $\overline{\cC}$ have irrational slope,
and in fact $\cC_+=\cC$ is an open convex cone.
To construct a $\Gamma_0$-invariant \rp\ decomposition $\cP$,
let $\Sigma$ be the convex hull of $\cC\cap\Phi(\aa)$.  The vertices
of $\Sigma$ form a countable set $\{v_j\}_{j\in\Z}$
which can be numbered so that the edges of $\Sigma$ are exactly
the line segments $\overline{v_jv_{j+1}}$.  If we let
$\sigma_j$ be the relative interior of the cone on $\overline{v_jv_{j+1}}$,
and let $\tau_j$ be the relative interior of the cone on
$v_j$, then $\cP:=\{\sigma_j\}_{j\in\Z}\cup\{\tau_j\}_{j\in\Z}$
is a $\Gamma_0$-invariant
\rp\ decomposition.  An explicit example of this
construction is illustrated on
p.~52 of \cite{AMRT}, reproduced as figure 1 of this paper.

The resulting partial compactification of $\cD/L$ adds a point $p_j$
for each $\sigma_j$, and a curve $B_j\cong\P^1$ for each $\tau_j$,
with $B_j\cap B_{j+1}=p_j$.  This can be pictured as an
``infinite chain'' of $\P^1$'s, as in the top of figure 2 (which
is also reproduced from \cite{AMRT}, p.~46).
The generator $[\diag(\varepsilon,1)]$ of $\Gamma_0=\Gamma_{\text{par}}/\aa$
acts by sending $v_j$ to $v_{j+m}$ for some fixed $m$.  Taking the
quotient by $\Gamma_0$ leaves us with a ``cycle'' of rational curves,
of length $m$ (as depicted in the bottom of figure 2).
We arrive at Hirzebruch's description of the resolution of the cusps.

\begin{figure}[t]
\begin{center}
\begin{picture}(10,5)(0,0)
\put(2,3){\mbox{Insert here the figure from p.~46 of \cite{AMRT}}}
\end{picture}
\end{center}
\caption{}
\end{figure}

\bigskip

Conversely, suppose we are given a normal surface singularity
$p\in \overline{S}$
(with $\overline{S}$ a small neighborhood of $p$) which has a resolution
of singularities $f:T\to \overline{S}$ such that
$B:=f^{-1}(p)$ is a cycle of rational curves,
that is, $B=B_1+\cdots+B_m$ is
a divisor with normal crossings such that $B_j$ only meets
$B_{j\pm1}$, with subscripts calculated mod $m$.
Much of the structure above can be recovered from this information alone.
In fact, by a theorem of Laufer \cite{laufer} these singularities are
{\em taut}, which means that the isomorphism type is determined by
the resolution data.  We will work  out in detail some aspects of this
tautness,
in preparation for a general construction in the
next section.

The starting point is
Wagreich's calculation \cite{wagreich} of
the local fundamental group $\pi_1(\overline{S}-p)$ for such singularities,
which goes as follows.  Let $S:=\overline{S}-p=T-B$.  The natural
map
$\iota:\pi_1(S)\to\pi_1(T)$
induced by the inclusion $S\subset T$
is surjective.  Since $T$ retracts onto a cycle of $\P^1$'s, the
group $\pi_1(T)\cong\pi_1(S^1)$ is infinite cyclic, and the
universal cover $\widehat{T}$ of $T$ contains an infinite chain
$\widehat{B}=\cdots+\widehat{B}_j+\widehat{B}_{j+1}+\cdots$ of $\P^1$'s
lying over the cycle $B$.  The kernel of $\iota$ is
$\pi_1(\widehat{T}-\widehat{B})$, and by a result of Mumford \cite{mumford}
this is a free abelian group generated by loops around any pair of
adjacent components $\widehat{B}_j$, $\widehat{B}_{j+1}$ of $\widehat{B}$.

In this way, we recover the two steps of the quotient construction,
and the compactification $\widehat{T}$ of the intermediate quotient
$\widehat{T}-\widehat{B}$.
Let $\widehat{S}$ be the universal cover of $S$ (and of
$\widehat{T}-\widehat{B}$).
To complete the discussion of tautness, we should exhibit an
isomorphism between
$\widehat{S}$ and an open subset of $\h\times\h$, which descends to a
$\pi_1(T)$-equivariant map
$(\widehat{T}-\widehat{B})\to(\h\times\h)/\aa$.  The easiest way to do this is
to
consider an extra piece of structure on $p\in \overline{S}$:
a flat connection
on the holomorphic cotangent bundle $\Omega^1_{S}$.  We discuss
this structure, and how to use it to determine the
mapping from $\widehat{S}$ to $\h\times\h=\cD$, in the next section.
(To give a complete proof of Laufer's tautness result along these
lines, we would also need to show how the connection is to be constructed;
we will not attempt to do that here.)

\section{The flat connection} \label{flat}

Let $(L_{\Q},\cC,\Gamma_0)$ be
 an admissible triple,
with associated tube domain $\cD=L_{\R}+i\,\cC$
and discrete group
$\Gamma=L\semidirect\Gamma_0\subset\Aff(L)$. We will define a
flat connection on the holomorphic cotangent bundle of the
quotient space $\cD/\Gamma$.

The intermediate quotient space $\cD/L$ is an open subset of the
algebraic torus $L_{\C}/L=L\otimes_{\Z}\C^*\cong(\C^*)^{\rk(L)}$.
We identify the dual of the Lie algebra $\Lie(L_{\C}/L)^*$ of that
torus with the space of right-invariant one-forms on the group $L_{\C}/L$.
Any basis of $\Lie(L_{\C}/L)^*$, when regarded as a subset
 of the space of global sections of the sheaf
$\Omega^1_{L_{\C}/L}$,
freely generates that sheaf at any point.
We can therefore define a connection $\nabla_{\text{toric}}$ on
$\Omega^1_{L_{\C}/L}$, the {\em toric connection}, by the requirement that
$\nabla_{\text{toric}}(\alpha)=0$ for every $\alpha\in\Lie(L_{\C}/L)^*$.
Since the group $L_{\C}/L$ is abelian, the connection $\nabla_{\text{toric}}$
is flat.

The action of $\Aff(L)$ on $L_{\C}$ descends to an action of $\GL(L)$
on $L_{\C}/L$ which preserves the space of right-invariant one-forms.
In particular, the $\GL(L)$-action will be compatible with the toric
connection.  Thus, if we restrict $\nabla_{\text{toric}}$ to $\cD/L$,
it commutes with the action of $\Gamma_0$ and induces a connection on
the holomorphic cotangent bundle of
$(\cD/L)/\Gamma_0=\cD/\Gamma$, still denoted by $\nabla_{\text{toric}}$.

Let $\sigma\subset L_{\R}$ be the relative interior of a rational
polyhedral cone which is generated by a basis $\ell^1,\dots,\ell^r$
of $L$, and let $z_1,\dots,z_r$ be the
coordinates on $L_{\C}$ dual to $\{\ell^j\}$.  The one-forms
$\dlog w_j:=2\pi i\,dz_j$ are right-invariant one-forms on $L_{\C}/L$
which serve as a basis of $\Lie(L_{\C}/L)^*$.  If we compactify
the open set $\cD_\sigma/L\subset L_{\C}/L$ to $U_\sigma:=(\cD_\sigma/L)^-$,
then the forms $\dlog w_j$ extend to meromorphic one-forms on $U_\sigma$
with poles along the boundary $B_\sigma:=(\cD_\sigma/L)^--(\cD_\sigma/L)$.
In fact, the forms
$\dlog w_1,\dots,\dlog w_r$ freely generate
the sheaf $\Omega^1_{U_\sigma}(\log B_\sigma)$
as an $\cO_{U_\sigma}$-module.  The flat connection $\nabla_{\text{toric}}$
therefore extends
to a flat connection on $\Omega^1_{U_\sigma}(\log B_\sigma)$ for which
the $\dlog w_j$ are flat sections.
Note that the connection does {\em not}\/ acquire singularities along
the boundary, but extends as a regular connection to the sheaf of
logarithmic differentials.

If $\cP$ is a \rp\ decomposition of $\cC_+$, we get in this way
an extension of the
flat connection $\nabla_{\text{toric}}$ from
$\Omega^1_{\cD/L}$ to the sheaf of logarithmic differentials on
$\widehat{\cD}(\cP)/L$ with poles on the boundary
$(\widehat{\cD}(\cP)/L)-(\cD/L)$.  As this extended connection still
commutes with $\Gamma_0$, there is an induced extension of
$\nabla_{\text{toric}}$ from $\Omega^1_{\cD/\Gamma}$ to
$\Omega^1_{\widehat{\cD}(\cP)/\Gamma}(\log\cB)$, where
$\cB:= (\widehat{\cD}(\cP)/\Gamma)-(\cD/\Gamma)$.  This holds for
any Mumford-type semi-toric compactification.

The existence of this toric connection on $\cD/\Gamma$
depends in an essential way on $\Gamma$
being a group of affine-linear transformations of $L$.  If $\cD$
admits an action by
a larger group $\Gamma_{\text{big}}$ which includes discrete
symmetries that do not lie in $\Aff(L)$, then $\nabla_{\text{toric}}$
may fail
to descend to the quotient $\cD/\Gamma_{\text{big}}$.  For example,
if $L=\Z$ acts on the upper half plane $\h$ by translations, then
the associated flat connection $\nabla_{\text{toric}}$ has the
property that
$\nabla_{\text{toric}}(d\tau)=0$, where $\tau$ is the standard
coordinate on $\h$.  The flat section $d\tau$ is
invariant under translations
$\tau\mapsto\tau+n$, but if we
apply $\nabla_{\text{toric}}$ to the pullback of the flat section $d\tau$
under
the inversion
 $\tau\mapsto{-1}/\tau$
we get
\[\nabla_{\text{toric}}(\tau^{-2}\,d\tau)=-2\tau^{-3}\,d\tau\otimes d\tau,\]
which is not $0$.  In particular, the
 connection $\nabla_{\text{toric}}$ does not descend to
the $j$-line $\h/\SL(2,\Z)$.

\bigskip

We now want to
explain how the abstract knowledge of the flat connection
$\nabla_{\text{toric}}$
and of a Mumford-type semi-toric compactification of $\cD/\Gamma$
can be used to recover the structure of $\cD$ and of $\Gamma$.
Suppose we are given a complex manifold $T$,
a divisor with normal crossings $B$ on $T$, and a flat connection
$\nabla$ on $\Omega^1_T(\log B)$.
By the usual equivalence between
flat connections and local systems \cite{RSP}, the flat sections
of $\nabla$ determine a local system $\E$ on T.  Such a local
system is specified by giving its fiber $E$ at a fixed base point $\star$
(which we choose to lie in $T-B$),
together with a representation of $\pi_1(T,\star)$ in $\GL(E)$.

We first restrict the connection and the local system to $T-B$.
If we pass to the universal cover $\widehat{S}$ of $T-B$, the
flat sections give a global trivialization of the bundle
$\E\otimes\cO_{\widehat{S}}=\Omega^1_{\widehat{S}}$.  There is a
natural map $\operatorname{int}_\star:\widehat{S}\to E^*$
which sends $s\in\widehat{S}$
to the functional
\[\alpha\mapsto\int_\star^s\widehat{\alpha},\]
where $\widehat{\alpha}$ is the unique flat section of $\E$ (a
holomorphic 1-form on $\widehat{S}$)
such that $\widehat{\alpha}|_\star=\alpha\in E$.
(Notice that if we vary the basepoint $\star$, we simply shift the
image of the map by some constant vector in $E^*$.)

On the other hand, if we consider $\nabla$ on $T$ and
 pass to the universal cover $\widehat{T}$ of $T$, the flat
sections of $\E$ will trivialize the bundle
$\Omega^1_{\widehat{T}}(\log\widehat{B})$, where
 $\widehat{B}$
is a divisor with normal crossings in $\widehat{T}$,
the inverse image of $B\subset T$.
We once again encounter the
intermediate quotient space $\widehat{T}-\widehat{B}$, and its
partial compactification $\widehat{T}$.

At any maximal-depth normal crossing point $p$
of $\widehat{B}\subset\widehat{T}$, let $v_{j}=0$
define the $j^{\text{th}}$ local component $B_j$ of the boundary at $p$.
There
is a unique flat section $\widehat{\alpha}_j$ of
$\Omega^1_{\widehat{T}}(\log\widehat{B})$,
defined locally near $p$, such that $\widehat{\alpha}_j-\dlog v_j$
vanishes at $p$.  It follows that $\widehat{\alpha}_j-\dlog v_j$
is a holomorphic one-form in a neighborhood of $p$, and so that
$\widehat{\alpha}_1,\dots,\widehat{\alpha}_r$ is a basis for
(flat) local sections of
$\Omega^1_{\widehat{T}}(\log\widehat{B})$.  Using the global trivialization,
we may regard each $\alpha_j:=\widehat{\alpha}_j|_\star$
 as an element of $E$.
We let $L_p\subset E^*$  be the lattice spanned by the dual basis
$\ell^1,\dots,\ell^r$ to $\alpha_1,\dots,\alpha_r$, and let
$\sigma_p\subset L_p\otimes\R$ be the relative interior of the cone
generated by $\ell^1,\dots,\ell^r$.

If we are to recover the structure of the semi-toric compactification,
we need a certain compatibility among the $L_p$'s and the $\sigma_p$'s:
they should be related to a common lattice and a common cone, independent
of $p$.  We formalize this as follows.

\begin{definition} \label{compat}
We call $(T,B,\nabla)$\/ {\em compatible} provided that
\begin{enumerate}
\item
each component of $B$ contains at least one maximal-depth normal crossing
point,

\item
the lattices $L_p$ for maximal-depth normal crossing points $p$
all coincide with a common lattice $L\subset E^*$,

\item
the natural map $\operatorname{int}_\star:\widehat{S}\to E^*=L_{\C}$
descends to a map
$(\widehat{T}-\widehat{B})\to(L_{\C}/L)$ which induces an isomorphism
of fundamental groups, and

\item
the collection $\cP$ of relative interiors of faces of the $\sigma_p$'s
is a \lrp\ decomposition of a strongly convex cone $\cC_+$.
\end{enumerate}
\end{definition}

Suppose that $(T,B,\nabla)$ is compatible,  let $\cC$ be the
interior of $\cC_+$, and let $\cD=L_{\R}+i\,\cC$.  The action
of $\pi_1(T)$ on $L_{\C}$ permutes the set of maximal-depth
normal crossing points of $B\subset T$,
and so preserves $\cP$ and $\cC$.  Thus, $\Gamma:=\pi_1(T-B)$
acts on $\cD$, and there is an induced map $(T-B)\to(\cD/\Gamma)$.

We can now recover the compactification $T$ from this data
(or at least its structure in codimension one).
For any maximal-depth
normal crossing boundary point $p$ of $\widehat{B}\subset\widehat{T}$,
there is a neighborhood $U_p$ of $p$ in $\widehat{T}$ and a natural
extension of the induced map
$U_p\cap(\widehat{T}-\widehat{B})\to L_{\C}/L$ to
a map $U_p\to\widehat{\cD}(\cP)/L$.  We cannot tell from the behavior
of these extensions what happens at
``interior'' points of boundary components (those which do not lie
in any $U_p$), but we {\em can}\/
conclude that there is a meromorphic map
$\widehat{T}\to\widehat{\cD}(\cP)/L$ which does not blow down any
boundary components.  This map is $\pi_1(T)$-equivariant, so it
descends to a map $T\to\widehat{\cD}(\cP)/\Gamma$.

\section{Moduli spaces of sigma-models} \label{sigma}

A {\em Calabi-Yau manifold}\/ is a compact connected orientable
 manifold $X$ of dimension $2n$
which admits Riemannian metrics whose holonomy is contained in
$\SU(n)$.\footnote{There is some confusion in the literature about
whether ``Calabi-Yau'' should mean that the holonomy is precisely
$\SU(n)$, or simply contained in $\SU(n)$.  In this paper, we adopt
the latter interpretation.}
Given such a metric, there exist complex structures on $X$ for
which the metric is K\"ahler.
The holonomy condition
is equivalent to requiring that this K\"ahler metric be Ricci-flat
(cf.~\cite{beauville}).
On the other hand, if we are given a
complex structure on a Calabi-Yau manifold, then by the theorems of
Calabi \cite{calabi} and Yau \cite{yau}, for each K\"ahler
metric $\widetilde{g}$ there is a unique
Ricci-flat K\"ahler metric $g$ whose
K\"ahler form is in the same de Rham cohomology class as that of
$\widetilde{g}$.
(We have implicitly used the topological consequence of Ricci-flatness:
Calabi-Yau manifolds  have vanishing first Chern class.)

Examples of Calabi-Yau manifolds are provided by the differentiable
manifolds underlying smooth complex projective varieties
 with trivial canonical bundle.  One can apply Yau's theorem to
a K\"ahler metric coming from a projective embedding in order to produce
a metric with holonomy contained in $\SU(n)$,
where $n$ is the complex dimension of the variety.  As explained in
\cite{beauville}, if the Hodge numbers $h^{p,0}$ vanish for $0<p<n$, then
the holonomy of this metric is precisely $\SU(n)$.

Physicists have constructed a class of conformal field theories called
\NLSM s on Calabi-Yau manifolds $X$ (cf.~\cite{GSW,hubsch}).
We consider here an approximation to those
theories, which should be called
``one-loop semiclassical \NLSM s''.
Such an object is determined by the data of a Riemannian metric
$g$ on $X$ whose holonomy is contained in $\SU(n)$ together with
the de~Rham cohomology class $[b]\in H^2(X,\R)$ of a real closed
$2$-form $b$ on $X$.

Two such pairs
$(g,b)$ and $(g',b')$ will
determine isomorphic conformal field theories if there is a
diffeomorphism $\varphi:X\to X$ such that $\varphi^*(g')=g$,
and $\varphi^*([b'])-[b]\in H^2(X,\Z)$.
It is therefore natural to regard
the class of $[b]$ in
$H^2(X,\R)/H^2(X,\Z)$
as the fundamental datum.
We denote this class by
$[b]\bmod\Z $.

The set of all isomorphism classes of such pairs we call the
{\em one-loop semiclassical \NLSM\ moduli space}, or simply the
{\em sigma-model moduli space}\/  (for short).
This may differ from the actual {\em conformal field theory moduli space},
both because
there may be additional
isomorphisms of conformal field theories which are not visible in this
geometric interpretation,
and also because there may be deformations of the \NLSM\ as a conformal
field theory which do not have a sigma-model interpretation on $X$
(cf.~\cite{mmm,phases}).
For our present purposes, we ignore
these more delicate questions about the
 conformal field theory moduli space, and
 concentrate on the sigma-model moduli space we have defined
above.

We focus attention in this paper on the case in which the holonomy
of the metric $g$ is precisely $SU(n)$, $n\ne2$.  For each such metric,
there are exactly two complex structures
on $X$ for which the metric is
K\"ahler (complex conjugates of each other).\footnote{More generally,
as we will show elsewhere,
if $h^{2,0}(X)=0$ there are only a finite number of complex structures
for which $g$ is K\"ahler.  The number depends on the decomposition
of the holonomy representation into irreducible pieces.}
Thus, there is
 a natural map from a double cover of
the sigma-model moduli space to
the usual ``complex structure moduli space'', given by assigning to
$(g,b)$ one of the two complex structures for which $g$ is K\"ahler.  The
fibers of this map can be described as follows.  If we fix a complex
structure on $X$, then the corresponding fiber consists of all
$B+i\,J\bmod\Z\in H^2(X,\C)/H^2(X,\Z)$
(modulo diffeomorphism)
with $B$ denoting the class $[b]$, for which $J$ is
the cohomology
class of a K\"ahler form.  (The metric $g$ is  uniquely determined
by $J$, by Calabi's theorem.)
This quantity
$B+i\,J\bmod\Z$ describes the ``extra structure'' $S$ which was alluded to
in \cite{guide}.
This is often called the {\em complexified K\"ahler structure}\/
on $X$ determined by $(g,b)$.

The natural double cover of the sigma-model moduli space will be
locally a product near $(g,b)$, with the variations of complex structure
and of complexified K\"ahler structure describing the factors in
the product, provided that neither the K\"ahler cone
nor the group of holomorphic automorphisms
  ``jumps''
when the complex structure varies.  (The non-jumping
of the K\"ahler cone was shown to hold
by Wilson \cite{wilson} in the case of holonomy $\SU(3)$, when the complex
structure is generic.)
We will tacitly assume this local product structure, and separately
study the parameter spaces for the variations of complexified K\"ahler
structure and of complex structure.

With a fixed complex structure on $X$, the parameter space for
complexified K\"ahler structures on $X$
can be described
in terms of the K\"ahler cone $\cK$ of $X$, and the lattice
$L=H^2(X,\Z)/(\text{torsion})$.
We must  identify any pair of complexified K\"ahler structures
which differ by a diffeomorphism that fixes the complex structure,
that is, by an element of
the group $\Gamma_0=\Aut(X)$
of holomorphic automorphisms.
The natural parameter space for  pairs
$(g,b)$ such that $g$ is K\"ahler for the given complex structure
thus has the form $\cD/\Gamma$,
where $\cD=\{B+i\,J\suchthat J\in\cK\}$ and $\Gamma=L\semidirect\Gamma_0$
is the extension of $\Gamma_0$ by the lattice translations.
  This is exactly
the kind of space encountered in the first part of this paper: a tube domain
modulo a discrete symmetry group of affine-linear transformations
which includes a lattice acting by translations.

A common technique in the physics literature is to consider
what happens
along paths $\{tz\bmod\Gamma\}_{t\to\infty}$, which go
from $z\in\cD$ out towards infinity
in the tube domain.  Many aspects
of the conformal field theory can  be analyzed perturbatively
in $t$ along such paths.  It seems
reasonable to hope that such limits can be
described in a common framework, based on a single partial compactification
of $\cD/\Gamma$.
This hope (together with a bit of evidence, discussed below)
leads us to conjecture  that
$(L_{\Q},\cK,\Aut(X))$
is an admissible triple, in order that Looijenga's methods
could be applied to construct compactifications of $\cD/\Gamma$.
We formulate this conjecture more explicitly as follows.

\begin{ConeConjecture}
Let $X$ be a Calabi-Yau manifold
on which a complex
structure has been chosen, and suppose that
$h^{2,0}(X)=0$.
Let $L:=H^{2}(X,\Z)/\text{torsion}$, let $\cK$ be the K\"ahler
cone of $X$, let $\cK_+$ be the convex hull of $\overline{\cK}\cap L_{\Q}$,
and let $\Aut(X)$ be the group of holomorphic automorphisms of $X$.
Then there exists a rational polyhedral cone $\Pi\subset\cK_+$
such that $\Aut(X).\Pi=\cK_+$.
\end{ConeConjecture}

The K\"ahler cone of $X$ can have a rather complicated structure, analyzed
in the case $n=3$
by Kawamata \cite{kawamata} and Wilson \cite{wilson}.
Away from classes of triple-self-intersection zero, the closed cone
$\overline{\cK}$ is locally rational polyhedral, but the rational
faces may accumulate towards points with vanishing triple-self-intersection.
The cone conjecture predicts that while the closed cone $\overline{\cK}$ of
$X$ may have infinitely many edges, there will only be finitely many
$\Aut(X)$-orbits of edges.  Other finiteness predictions which follow
from the cone conjecture include finiteness of the set of fiber space
structures on $X$, modulo automorphisms.

Many of the large classes of examples, such as toric hypersurfaces, have
K\"ahler cones $\cK$ such that $\cK_+=\overline{\cK}$ is a rational polyhedral
cone.  For these, the cone conjecture automatically holds.
A nontrivial case of the cone conjecture---Calabi-Yau threefolds which are
fiber products of generic rational elliptic surfaces with section (as
studied by Schoen \cite{schoen})---has been checked by Grassi
and the author \cite{grassi-morrison}.
In addition, Borcea \cite{borcea} has verified the
finiteness of $\Aut(X)$-orbits of edges of $\overline{\cK}$ in another
nontrivial example, and
Oguiso \cite{oguiso} has discussed finiteness of $\Aut(X)$-orbits of
fiber space structures in yet another example.\footnote{Neither of these
constitutes a complete verification of the cone conjecture for the
threefold in question.}
All three examples involve cones with an
infinite number of edges.

For any $X$ for which the cone conjecture holds,
the K\"ahler parameter space $\cD/\Gamma$
 will admit both a \SBB-type
``minimal'' compactification, and smooth compactifications
of Mumford type built out
of many cones $\sigma\subset\cK$ as above.

\section{Additional structures on the moduli spaces} \label{additional}

Of particular interest to the physicists
studying \NLSM s
has been the ``large radius
limit'' in the K\"ahler parameter space.  This is typically analyzed in
the physics literature as follows (cf.~\cite{tsm,topgrav}).
The quantities of physical interest
will be
invariant under translation by $L$.  Many such quantities
vary holomorphically
with parameters, and their Fourier expansions take the form
\renewcommand{\theequation}{*}
\begin{equation} \label{eq}
\sum_{\eta\in L^*}c_\eta \,e^{2\pi i\,z\cdot \eta}.
\end{equation}
The coefficients $c_\eta$ for $\eta\ne0$
are called {\em instanton contributions}\/
to the quantity \eqref{eq}, and
in many cases they can be
given a geometric interpretation which shows that they
vanish unless
$\eta$ is the class of an effective curve on $X$.  A
``large radius limit'' should be a point at which instanton
contributions to quantities like \eqref{eq} are suppressed
\cite{greene-plesser,al2}.

If we pick a basis $\ell^1,\dots,\ell^r$ of $L$ consisting of
vectors which lie in the closure of the K\"ahler cone,
write $\eta=\sum \eta^j\ell_j$
in terms of the basis $\{\ell_j\}$ of $L^*$ dual to $\{\ell^j\}$,
and express
\eqref{eq} as a power series in $w_j:=\exp(2\pi i\,z_j)$,
where $\{z_j\}$ are coordinates dual to $\{\ell^j\}$,
then the series expansion
\renewcommand{\theequation}{**}
\begin{equation} \label{eqq}
\sum_{\eta\in L^*}c_\eta \,w_1^{\eta^1}\cdots w_r^{\eta^r}
\end{equation}
involves only terms with nonnegative exponents \cite{asp-mor}.
If convergent,\footnote{From a rigorous mathematical point of view,
the Fourier coefficients
$c_\eta$ can often
be defined and calculated, but no convergence properties
of the series \eqref{eq} or \eqref{eqq} are known.}
 this will define a function on $(\cD_\sigma/L)^-$,
where $\sigma$ is the relative interior of the cone generated by
$\ell^1,\dots,\ell^r$.
Thus, approaching the distinguished limit point of $\cD_\sigma/L$
(where all $w_j$'s approach $0$) suppresses the
instanton contributions, so the distinguished limit point is
a good candidate for the large radius limit.
We can repeat this construction for any cone $\sigma\subset\cK$
which is the relative interior of a cone
generated by a basis of $L$, obtaining partial
compactifications which include large radius limit points
for paths that lie in various cones $\sigma$.

Among the ``quantities of physical interest'' to which this analysis
is applied are a collection
of multilinear
maps of cohomology groups
called {\em three-point functions}.  These maps should depend
on the data $(g,b)$,
and should vary holomorphically with both complex structure and complexified
K\"ahler structure parameters.
Certain of these three-point functions
(related to Witten's ``$A$-model'' \cite{witten}) would depend only
on the complexified K\"ahler structure,
while others (related to Witten's ``$B$-model'') would
depend only on the complex structure.  The $B$-model three-point
functions can be mathematically interpreted in terms of the variation
of Hodge structure, or period mapping,
induced by varying the complex structure on
the Calabi-Yau manifold
\cite{cecotti2,guide,gmp}.

In \cite{htams}, we discuss a mathematical version
of the $A$-model three-point functions, expressed as formal power series
near the distinguished limit  point associated to
the relative interior $\sigma$ of a rational polyhedral
 cone  generated by a basis of $L$.
(The coefficients $c_\eta$ of this power series are
derived from the numbers of rational curves on $X$ of various degrees.)
The choice of $\sigma$ is an additional piece of data
in the construction which we call
a {\em framing}.

These formal power series representations of $A$-model three-point functions
can be regarded as defining a formal
degenerating variation of Hodge structure, which we call the
{\em framed $A$-variation of Hodge structure}\/  with framing $\sigma$.
Now there are manipulations of these formal series which suggest that
the underlying convergent three-point functions (if they exist) will
not depend on the choice of $\sigma$ and will be invariant under the
action of $\Aut(X)$.

We must refer the reader to \cite{htams} for the precise definition of
framed $A$-variation of Hodge structure.  But for reference, we would like
to state here a conjecture which suggests how the various framed
$A$-variations of Hodge structure will fit together, along the
lines being discussed in this paper.

\begin{ConvergenceConjecture}
Suppose that $X$ is a Calabi-Yau manifold with $h^{2,0}(X)=0$,
endowed with a complex
structure, which
satisfies the cone conjecture.  Let $L:=H^2(X,\Z)/\text{torsion}$,
let $\cK$ be the K\"ahler cone of $X$, let $\cD:=L_{\R}+i\,\cK$
be the associated tube domain, and let
$\Gamma:=L\semidirect\Aut(X)$.
Then there is a neighborhood
$U$ of the $0$-dimensional stratum $\widehat{\cD}(\cK)$ in the \SBB-type
compactification
$\widehat{\cD}(\cP_{\text{SBB}})/\Gamma$,
and a variation
of Hodge structure on $U\cap(\cD/\Gamma)$, such that for any
$\sigma\subset\cK$ which is the relative interior of a
\rp\ cone
$\overline{\sigma}\subset\cK_+$  generated by a basis of $L$,
the induced formal degenerating
variation of Hodge structure at the distinguished limit point of
$\cD_{\sigma}/L$ agrees with the
framed $A$-variation of Hodge structure
with framing $\sigma$.
\end{ConvergenceConjecture}

If this variation of Hodge structure exists, we call it the
{\em $A$-variation of Hodge structure}\/ associated to $X$.

\section{Maximally unipotent boundary points} \label{MUBP}

In the previous section, we discussed how to let the
complexified K\"ahler
parameter $B+i\,J$ approach infinity, analyzing certain partial
compactifications and boundary points of the sigma-model moduli
space in the $B+i\,J$ directions.  We now turn to compactifications
and boundary points in the transverse directions---the directions
obtained by varying the complex structure on the Calabi-Yau manifold.
We consider what
happens when the complex structure degenerates.

The local moduli spaces of complex structures on Calabi-Yau manifolds
are particularly well-behaved, thanks to a theorem
of Bogomolov \cite{bogomolov}, Tian \cite{tian}, and Todorov \cite{todorov},
which guarantees that all first-order deformations are unobstructed.
In particular, there will be a local family of deformations of a given complex
structure
for which the Kodaira-Spencer map is an isomorphism.
More generally, we consider arbitrary
 families $\pi:\cY\to S$ of complex structures on a fixed Calabi-Yau
manifold $Y$, by which we mean:  $\pi$ is a proper and smooth
map between connected
complex manifolds, and all fibers $Y_s:=\pi^{-1}(s)$ are diffeomorphic
to $Y$.  We will often assume that the Kodaira-Spencer map
is an isomorphism at every point $s\in S$, so that $S$ provides
 good local moduli spaces for the fibers $Y_s$.

To study the behavior when the complex structure degenerates,
 we partially compactify the parameter space $S$ to $\overline{S}$.
There is a class of boundary points on $\overline{S}$
of particular interest from
the perspective of conformal field theory.
According to the interpretation of \cite{guide,htams},
these points can be
identified by the monodromy properties of the associated variation
of Hodge
structure\footnote{The variation of Hodge structure in question is the
usual geometric one (cf.~\cite{transcendental})
associated to a variation of complex structure.
These might be called ``$B$-variations of Hodge structure'' by analogy
with the previous section.}
near $p\in\overline{S}$.
We first review from
\cite{htams} these monodromy properties for normal crossing
boundary points, and then extend the definition
to a wider class of compactifications and boundary points.

Let $p$ be a maximal-depth normal crossing point of $B\subset\overline{S}$,
where $B:=\overline{S}-S$ is the boundary, assumed for the moment to
be a divisor with normal crossings.
Let $U$ be a small neighborhood of $p$
in $\overline{S}$, and write $B\cap U$ in the form
$B_1+\cdots+ B_r$.
If we fix a point $s\in U{\stminus}B$,
then each
local divisor $B_j$ gives rise to an
 monodromy transformation $T^{(j)}:H^n(Y_s,\Q)\to H^n(Y_s,\Q)$, which
is guaranteed to be quasi-unipotent by the monodromy theorem
\cite{monodromy}.

\begin{definition}
A maximal-depth normal crossing  point $p$ of $B\in\overline{S}$
is called a\/ {\em maximally unipotent point}\footnote{When
 $\dim(\cM)=1$,
this definition is equivalent to
the one given in \cite{guide}.} under the following conditions.
\begin{enumerate}
\item
The monodromy transformations $T^{(j)}$ around local boundary components
$B_j$ near $p$ are all unipotent.
\item
Let $N^{(j)}:=\log T^{(j)}$, let $N :=\sum a_jN^{(j)}$ for some $a_j>0$,
and define
\begin{eqnarray*}
W_0&:=&\Im(N ^n)\\
W_1&:=&\Im(N ^{n-1})\cap\Ker N \\
W_2&:=&\Im(N ^{n-2})\cap\Ker(N ^2).
\end{eqnarray*}
Then $\dim W_0=\dim W_1=1$ and $\dim W_2=1+\dim(\cM)$.
\item
Let $g^0,g^1,\dots,g^r$ be a basis of $W_2$ such that $g^0$ spans $W_0$,
and define $m^{jk}$ by
$N^{(j)}g^k=m^{jk}g^0$ for $1\le j,k\le r$.
Then $m:=(m^{jk})$ is an invertible matrix.
\end{enumerate}
(The spaces $W_0$ and $W_2$ are independent of the choice of coefficients
$\{a_j\}$
{\rm\cite{CK,deligne}},
and the invertibility of $m$ is independent of the choice
of basis $\{g^k\}$.)
\end{definition}

Given a maximally unipotent point $p\in\overline{S}$, we define the
{\em canonical logarithmic one-forms}\/
$\dlog q_j\in\Gamma(U,
\Omega_{\overline{S}}^1(\log B))$
at $p$ by
\[\frac1{2\pi i}\dlog q_j:=
d\left(\frac{\sum_{k=1}^r\langle g^k|\omega\rangle m_{kj}}
{\langle g^0|\omega\rangle}\right)\]
where $(m_{kj})$ is the inverse matrix of $(m^{jk})$, and $\omega$ is
a section of the sheaf $\Omega^n_{\cY/S}$ of
relative holomorphic $n$-forms on the family of complex structures
parameterized by $S$.
The elements $g^k\in H^n(Y_s,\Q)$ have been implicitly
extended to multi-valued sections of the local system $R^n\pi^*(\Q_{\cY})$
in order
to evaluate $\langle g^k|\omega\rangle$;
the monodromy measures the multi-valuedness of the resulting
(locally defined) holomorphic functions $\langle g^k | \omega \rangle$.
The fact that each $\dlog q_j$ as defined above
has a single-valued meromorphic extension
to $U$ follows from the nilpotent orbit theorem \cite{schmid}.
In \cite{htams} we show that the canonical
one-forms are independent of the choice of basis $\{g^k\}$,
and also of the choice of relative
$n$-form $\omega$; that for  any local defining
equation $v_j=0$ of $B_j$, the one-form $\dlog q_j-\dlog v_j$ extends
to a regular one-form on $U$; and that $\dlog q_1$, \dots, $\dlog q_r$
freely generate the locally free sheaf
$\Omega_{\overline{S}}^1(\log B)$
near $p$.

The canonical logarithmic one-forms can be integrated to produce
{\em quasi-canonical coordinates}\/
$q_1$, \dots, $q_r$ near $p$, but due to  constants
of integration, these coordinates are not unique.
That is, if we attempt to define
\[q_j=\exp\left(2\pi i\,
\frac{\sum_{k=1}^r\langle g^k|\omega\rangle m_{kj}}
{\langle g^0|\omega\rangle}\right)\]
we find that changing the basis $\{g^k\}$ will alter the $q_j$'s by
multiplicative constants (cf.~\cite{picard-fuchs}).
To specify
truly canonical coordinates, further conditions on the basis $\{g^k\}$
must be imposed, as discussed in \cite{guide,htams}.
For example, by demanding that $g^0$ span $W_0\cap H^n(Y_s,\Z)/\text{torsion}$
and that $g^0,\dots,g^r$ span $W_2\cap H^n(Y_s,\Z)/\text{torsion}$
we can reduce the ambiguity in the $q_j$'s to a finite number of choices.

With no ambiguity, we can
use the canonical logarithmic one-forms to produce a (canonical)
flat connection
$\flt$ on the holomorphic vector bundle
$\Omega_{U}^1(\log B)$
by declaring $\dlog q_1$, \dots, $\dlog q_r$ to be a basis for the
$\flt$-flat sections, that is, $\flt(\dlog q_j)=0$.  Notice that the
 connection $\flt$ is regular along the
boundary divisor $B$.
This connection is what we will use to extend the definition of maximally
unipotent to a more general case.

We now consider   partial compactifications
$\overline{S}$ of $S$ which are not necessarily smooth,
and whose boundary
is not necessarily a divisor with normal crossings.

\begin{definition} \label{def3}
Let $\Xi\subset\overline{S}{\stminus}S$ be a connected
subset of the boundary.
We say that $\Xi$ is\/ {\em maximally unipotent}  if there
is a neighborhood $V$ of  $\Xi$ in
$\overline{S}$
and a flat connection $\nabla_{\text{unip}}$ on
$\Omega^1_{V\cap S}$ such that for some
resolution of singularities $f:U\to V$ which is an isomorphism over $V\cap S$,
we have
\begin{enumerate}
\item
 the new boundary $B:=U{\stminus}f^{-1}(V\cap S)$ on $U$
is a divisor with normal
crossings,

\item
the flat connection $\nabla_{\text{unip}}$ extends to a connection on
$\Omega^1_U(\log B)$
(also denoted by $\nabla_{\text{unip}}$),

\item
for every maximal-depth normal crossing point $p$
of $B\subset U$, we have $\nabla_{\text{unip}}(\dlog q_j)=0$
for each
canonical logarithmic one-form $\dlog q_j$
at $p$,  and

\item
$(U,B,\nabla_{\text{unip}})$ is compatible in the sense of
definition \ref{compat}.
\end{enumerate}
We call $\nabla_{\text{unip}}$ the\/ {\em maximally unipotent connection
determined by $\Xi$}.
\end{definition}

Note that $\dlog q_1,\dots,\dlog q_r$ is a basis for the vector space
of local solutions of $\nabla_{\text{unip}}e=0$ near $p$.
By analytic continuation of solutions,
the connection $\nabla_{\text{unip}}$ is unique if it exists.
The requirement of compatibility is quite strong, essentially
guaranteeing that the structure of $\overline{S}$ near $\Xi$ resembles
that of a semi-toric compactification.

\section{Implications of mirror symmetry} \label{mirror}

Mirror symmetry \cite{greene-plesser}
predicts that Calabi-Yau manifolds should come in
pairs,\footnote{The most recent results \cite{mmm,phases} suggest that
it is {\em birational equivalence classes}\/ of Calabi-Yau manifolds which
come in pairs.}
with the r\^oles of variation of complex structure and of complexified
K\"ahler structure
being reversed between mirror partners.
We make a precise mathematical conjecture about mirror symmetry in
\cite{htams}, which can be stated as follows.

\begin{mirrorconjecture} {\rm (Normal Crossings Case)}
Let $Y$ be a Calabi-Yau manifold with $h^{2,0}(Y)=0$, and let
 $\pi:\cY\to S$ be a family of complex structures on $Y$ such that
the Kodaira-Spencer map is an isomorphism at every point.  Let
$S\subset\overline{S}$ be a partial compactification whose boundary
is a divisor with normal crossings.
To each maximally unipotent normal crossing boundary point $p$ in
$\overline{S}$
there is associated the following:
\begin{enumerate}
\item  a Calabi-Yau
manifold $X$ with $h^{2,0}(X)=0$,

\item
a lattice
$L$ of finite
index\footnote{The reason for allowing such an $L$ rather than insisting
 $H^2(X,\mbox{\footbbbfont Z})/\text{torsion}$
itself is that our basic
defining condition
on the family $S$---that the Kodaira-Spencer map be an isomorphism
at every point---is invariant under finite unramified
base change.  So we must allow
finite unramified covers of the parameter spaces.}
 in $H^2(X,\Z)/\text{torsion}$,

\item
the relative interior $\sigma\subset H^2(X,\R)$ of
a \rp\ cone which is generated by a basis $\ell^1,\dots,\ell^r$ of
$L$, and

\item
a map $\mu$ from a neighborhood of $p$ in $\overline{S}$
to $((H^2(X,\R)+i\,\sigma)/L)^-$, determined up to constants of integration
by the requirement that $\mu^*(\dlog w_j)$ is the canonical logarithmic
one-form $\dlog q_j$ on $\overline{S}$ at $p$ (as defined in
section \ref{MUBP}), where $z_1,\dots,z_r$
are coordinates dual to $\ell^1,\dots,\ell^r$, and $w_j:=\exp(2\pi i\,z_j)$,

\end{enumerate}
such that
\begin{itemize}
\item[a.]
$\sigma$ is contained in the K\"ahler cone for some complex structure
on $X$, and

\item[b.]
$\mu$ induces an isomorphism between the formally degenerating
geometric variation of Hodge structure at $p$ and the $A$-variation of
Hodge structure with framing $\sigma$ associated to $X$.
\end{itemize}
\end{mirrorconjecture}

Put more concretely, if we calculate the geometric variation of Hodge
structure near $p\in\overline{S}$
using appropriate quasi-canonical coordinates $q_j$, we should
produce power series expansions for $B$-model three-point functions (for $Y$)
whose coefficients agree with the  $c_\eta$ which
are derived from the numbers
of rational curves on $X$.  This is precisely the type of
calculation pioneered by Candelas, de la Ossa, Green, and Parkes \cite{pair}
in the case of the quintic threefold.

Given a family $\pi:\cY\to S$ of complex structures on $Y$, and
a partial compactification $\overline{S}$ of $S$,
if we move from point
to point along the boundary of $\overline{S}$, or if we
vary the compactification
$\overline{S}$ by blowing up the boundary, we can produce many maximally
unipotent normal crossing boundary points.
On the other hand, if $X$ is a mirror partner of $Y$ for which
 the cone and convergence conjectures hold,
there are many framed $A$-variations of Hodge structure
(with different framings) associated to $X$.
Given framings $\sigma$ and $\sigma'$ which belong to \rp\
decompositions $\cP$ and $\cP'$, respectively, there is always a common
refinement $\cP''$ of these decompositions.  Geometrically, the
corresponding compactification $\widehat{\cD}(\cP'')/\Gamma$ is a
blowup of both $\widehat{\cD}(\cP)/\Gamma$ and $\widehat{\cD}(\cP')/\Gamma$.
Analytic continuation on the common blowup $\widehat{\cD}(\cP'')/\Gamma$
from a point in the inverse image of
$\cD(\sigma)$ to one in the inverse image of $\cD(\sigma')$ will give
an isomorphism of the $A$-variations of Hodge structure.

The various maximally unipotent normal crossing boundary points will
 (conjecturally) lead to
many mirror isomorphisms.  We wish to fit
these various mirror isomorphisms together.  In fact, the mirror
symmetry isomorphism is expected by the physicists to extend
to an isomorphism between
the full conformal field theory moduli spaces, and so, presumably,
to compactifications as well. Thus, the structure of the semi-toric
compactifications which are natural from the point of view of
variation of complexified K\"ahler structure on $X$ should
be reflected in the structure of
compactifications of the complex structure moduli
space $\cM_Y$ of $Y$.

This philosophy suggests two things about the compactified parameter spaces
$\overline{S}$ of complex structures on $Y$.  First, there should
be a compatibility between compactification points whose mirror families
are associated to the {\em same}\/ space $X$, and the {\em same}\/ K\"ahler
cone $\cK$.  In fact, we should be able to extend our mathematical
mirror symmetry conjecture to arbitrary maximally unipotent subsets
of the boundary for {\em any}\/ compactification, not just ones
whose boundary is a divisor with normal
crossings.  And second, there should be some kind of
{\em minimal}\/ compactification of the coarse moduli space
$\cM_Y$ of complex structures on $Y$, whose mirror compactified
family would be the \SBB-type compactification of $\cD/\Gamma$.

The compatibility between compactifications can be recognized
 by means of the flat connection
$\nabla_{\text{unip}}$ which we used to identify maximally unipotent
subsets of the boundary.
We extend our mirror symmetry conjecture to the general case as follows.

\begin{mirrorconjecture} {\rm (General Case)}
Let $Y$ be a Calabi-Yau manifold with $h^{2,0}(Y)=0$,  let
 $\pi:\cY\to S$ be a family of complex structures on $Y$ such that
the Kodaira-Spencer map is an isomorphism at every point, and let
$S\subset\overline{S}$ be a partial compactification.
To each maximally unipotent connected subset $\Xi$ of the boundary
$\overline{S}{\stminus}S$
there is associated the following:
\begin{enumerate}
\item  a Calabi-Yau
manifold $X$ satisfying the cone and convergence conjectures,

\item
a subgroup $\Gamma\subset\Aff(H^2(X,\R))$ whose translation subgroup
$L$ is a lattice of finite index in $H^2(X,\Z)/\text{torsion}$,

\item
a \lrp\ decomposition $\cP$ of a cone $\cC_+$ (which coincides with the
convex hull of $\overline{\cC_+}\cap L_{\Q}$) that is invariant
under the group $\Gamma_0:=\Gamma/L$, and

\item
a map $\mu$ from a neighborhood $U$ of $\Xi$ in $\overline{S}$
to
$\widehat{\cD}(\cP)/\Gamma$,
 determined up to constants of integration
by the requirement that
the flat connection $\nabla_{\text{toric}}$ on $\cD/\Gamma$ pulls back
to  $\nabla_{\text{unip}}$ on $U\cap S$,
where $\nabla_{\text{unip}}$ is the maximally unipotent connection
determined by $\Xi$,
\end{enumerate}
such that
\begin{itemize}
\item[a.]
for some complex structure on $X$, the interior $\cC$ of $\cC_+$
is contained in the K\"ahler cone and $\Gamma_0$ is contained in the
group of holomorphic automorphisms, and

\item[b.]
$\mu$ induces an isomorphism between the
geometric variation of Hodge structure over $U\cap S$
 and the $A$-variation of
Hodge structure associated to $X$.
\end{itemize}
\end{mirrorconjecture}

{\em A priori}, the map $\mu$ determined by compatibility of the
connections would only be a meromorphic map; we are asserting that
it is in fact regular, and a local isomorphism.

There is one further refinement of this conjecture which could be
made: we could demand that the map $\mu$ also respect the quasi-canonical
coordinates determined by choosing integral bases $g^0,\dots,g^r$.
This would reduce the ambiguity in the choice of $\mu$ to a finite
number of choices, but would require a compatibility among such integral
quasi-canonical coordinates at various boundary points.

\bigskip

Finally, suppose that $\cM_Y$ is the coarse moduli space for complex
structures on a Calabi-Yau
variety  $Y$ such that $h^{2,0}(Y)$.
(This coarse moduli space is known to exist as a quasi-projective
variety, once we have specified a polarization, thanks to a theorem
of Viehweg \cite{viehweg}.)
In this case, we conjecture the existence of
a \SBB-style compactification, as follows.

\begin{minconjecture}
There is a partial compactification $(\overline{\cM_Y})_{\text{SBB}}$
of the coarse moduli space
$\cM_Y$ with distinguished boundary points $p_1,\dots,p_k$ which are
maximally unipotent, such that the data associated by the mathematical mirror
symmetry conjecture to $p_j$ consists of:
(1) a Calabi-Yau manifold $X_j$ (with a complex structure specified that
determines the group $\Aut(X_j)$ of holomorphic automorphisms and the
K\"ahler cone $\cK_j$ of $X_j$), (2) the group
\[\Gamma_j:=(H^2(X_j,\Z)/\text{torsion})\semidirect\Aut(X_j),\]
and (3) the \lrp\ decomposition $\cP_j$ which is the \SBB\ decomposition
$\cP_{\text{SBB}}$ of the cone $(\cK_j)_+$ (the convex hull of
$\overline{\cK}_j\cap H^2(X_j,\Q)$).
\end{minconjecture}

\section{Mumford cones and Mori cones} \label{mumford-mori}

In the fall of 1979, Mori  lectured at Harvard on his then-new
results \cite{mori-cone}
on the cone of effective curves.  In order to show that
his theorem about local finiteness of extremal rays fail when the
canonical bundle is numerically effective, he gave an example.
(A similar example appears in a Japanese expository paper he wrote
a few years later, which  has since been translated into English
\cite{mori-sugaku}.)
The example was of an abelian surface with real multiplication,
that is, one whose endomorphism algebra contains the ring of
integers $\O_K$ of a real quadratic field $K$.
For such a surface $X$, the N\'eron-Severi group
$L:=H^{1,1}(X)\cap H^2(X,\Z)$ is a lattice of rank $2$.
The K\"ahler cone of $X$
lies naturally in $L_{\R}$, and is an open cone $\cK$ bounded by
two rays whose slopes are irrational numbers in the field $K$
 (cf.~\cite{kuga-satake,vandergeer}).
Rays through classes of ample divisors $[D]\in L\cap\cK$ can be found
which are arbitrarily close to the boundary, but the boundary
is never reached.
This phenomenon indicated that Mori's results on the structure of the
dual cone $\cK^\vee\subset H_2(X,\R)$ could not be extended to the
case of abelian surfaces.
The picture Mori drew for this example was remarkably similar to
figure~1.

The Hilbert modular surfaces in fact serve as moduli spaces for
abelian surfaces with endomorphisms of this type
(cf.~\cite[Chap.~IX]{vandergeer}), although a bit more data must be
specified, which determines the group $\Gamma$.
Now Mumford's figure~1 was drawn in some auxiliary space
being used to describe this ``complex structure moduli space'', while
Mori's version of figure~1 depicted the K\"ahler cone in $H^{1,1}$,
and so is related to ``complexified K\"ahler moduli'' of the surfaces.
The setting
is not quite the same as the one in the present paper, since $h^{2,0}\ne0$.
However, mirror symmetry for complex tori does predict that
each cusp in
 the complex structure moduli space will be related to the K\"ahler
moduli space for the abelian varieties parametrized by some
$\h\times\h/\Gamma$, with the $\Gamma$ determined by the cusp.
(This is not completely
clear from the literature; I will return to this point in a subsequent
paper.) In fact,
under this association the Mumford cone from figure 1
corresponds precisely to the (dualized)
Mori cone.

Mirror symmetry might have been anticipated by mathematicians
had anyone noticed the striking similarity between these two pictures
back in 1979!

\section*{Acknowledgments}

I would like to thank P.~Aspinwall,
R.~Bryant,
P.~Deligne, G.~Faltings,
A.\ Grassi, B.~Greene, R.~Hain, S.~Katz, R.~Plesser, L.~Saper and E.~Witten
for useful discussions.  This work was partially supported by
NSF Grant  DMS-9103827 and by an American Mathematical Society
Centennial Fellowship.

\makeatletter \renewcommand{\@biblabel}[1]{\hfill#1.}\makeatother
\newcommand{\bysame}{\leavevmode\hbox to3em{\hrulefill}\,}

\trailer

\end{document}